\pdfoutput=1
\RequirePackage{ifpdf}
\ifpdf %We are running pdfTeX in pdf mode
\documentclass[pdftex]{sigma}
\else
\documentclass{sigma}
\fi

\newtheorem{thm}{Theorem}[section]
\newtheorem{lem}[thm]{Lemma}
\newtheorem{prop}[thm]{Proposition}
\newtheorem{cor}[thm]{Corollary}
\theoremstyle{definition}
\newtheorem{exam}[thm]{Example}
\newtheorem{defn}[thm]{Definition}

\numberwithin{equation}{section}

\newcommand\hyper[5]{
  {}_{#1}F_{#2} \left(
    \begin{matrix}
      #3\\
      #4\\
    \end{matrix}
    ; #5
    \right)
}
\newcommand\nth{n^{\mathrm{th}}}

\newcommand\mk[1]{\hat{#1}}
\newcommand\F[2]{{}_{#1}F_{#2}}

\begin{document}

\allowdisplaybreaks

\renewcommand{\thefootnote}{$\star$}

\renewcommand{\PaperNumber}{039}

\FirstPageHeading

\ShortArticleName{The Combinatorics of Associated Laguerre Polynomials}

\ArticleName{The Combinatorics\\
of Associated Laguerre Polynomials\footnote{This paper is a~contribution to the Special Issue on Exact Solvability and Symmetry
Avatars in honour of Luc Vinet.
The full collection is available at
\href{http://www.emis.de/journals/SIGMA/ESSA2014.html}{http://www.emis.de/journals/SIGMA/ESSA2014.html}}}

\Author{Jang Soo KIM~$^\dag$ and Dennis STANTON~$^\ddag$}

\AuthorNameForHeading{J.S.~Kim and D.~Stanton}

\Address{$^\dag$~Department of Mathematics, Sungkyunkwan University, Suwon 440-746, South Korea}
\EmailD{\href{jangsookim@skku.edu}{jangsookim@skku.edu}}

\Address{$^\ddag$~School of Mathematics, University of Minnesota, Minneapolis, MN 55455, USA}
\EmailD{\href{stanton@math.umn.edu}{stanton@math.umn.edu}}

\ArticleDates{Received January 30, 2015, in f\/inal form May 06, 2015; Published online May 11, 2015}

\Abstract{The explicit double sum for the associated Laguerre polynomials is derived combinatorially.
The moments are described using certain statistics on permutations and permutation tableaux.
Another derivation of the double sum is provided using only the moment generating function.}

\Keywords{associated Laguerre polynomial; moment of orthogonal polynomials, permutation tableau}

\Classification{05E35; 05A15}

\renewcommand{\thefootnote}{\arabic{footnote}}
\setcounter{footnote}{0}

\section{Introduction}%\label{sec:introduction}

The classical orthogonal polynomials of Hermite, Laguerre, and Jacobi have explicit representations as terminating hypergeometric
series.
In the Askey scheme of orthogonal polynomials, each polynomial is given by either a~hypergeometric or basic hypergeometric
series.
These polynomials have been studied extensively combinatorially.

The associated versions of the classical polynomials are def\/ined by replacing~$n$ by $n+c$ in the three term recurrence relation,
where~$c$ is an arbitrary parameter.
The associated classical polynomials may be expressed as terminating hypergeometric double sums.
For example, the coef\/f\/icient of $x^k$ in the associated Laguerre polynomial of degree~$n$ is given by a~terminating $\F32$
series.
The analogous expansion for the associated Askey--Wilson polynomials uses a~${}_{10}\phi_9$ series, see~\cite{IR}.

The combinatorics of the classical orthogonal polynomials extends easily to include the association parameter~$c$.
However it is not clear how these double sums occur in the explicit representations.
The purpose of this paper is to derive these representations in the Laguerre case, Theorems~\ref{doublesumthm}
and~\ref{doublesumthm2}, by giving combinatorial interpretations for the polynomials.
We also f\/ind combinatorial expressions for the moments, see Theorems~\ref{L2m} and~\ref{L1m}.
We use the standard methods of Motzkin paths and permutation statistics to carry this out.
Along the way we show that certain statistics on permutation tableaux describe the moments, Theorem~\ref{PT}.

As a~bonus, in Section~\ref{sec:proof-via-moment} we verify an explicit double sum expression for associated Laguerre polynomials
from knowledge of the moment generating function.
An explicit measure or a~verif\/ication of the 3-term recurrence is not necessary.
This method may be useful in other settings to f\/ind explicit representations for polynomials when we know neither the measure nor
the moments.

\section{The associated Laguerre polynomials}\label{sec:assoc-lagu-polyn}

The monic Laguerre polynomials $p_n(x)=(-1)^nn!L_n^\alpha(x)$, whose weight function is
\begin{gather*}
x^\alpha e^{-x}dx,
\qquad
x\ge 0,
\end{gather*}
may be def\/ined by the three term recurrence relation
\begin{gather}
\label{3rr}
p_{n+1}(x)=(x-b_n)p_n(x)-\lambda_n p_{n-1}(x),
\qquad
p_0(x)=1,
\qquad
p_{-1}(x)=0,
\end{gather}
where
\begin{gather*}
b_n=2n+\alpha+1,
\qquad
\lambda_n=n(n+\alpha).
\end{gather*}

The associated Laguerre polynomials satisfy the same recurrence with~$n$ replaced by $n+c$,
\begin{gather*}
b_n=2n+\alpha+2c+1,
\qquad
\lambda_n=(n+c)(n+\alpha+c).
\end{gather*}
For the purposes of this paper it is convenient to replace the parameters~$\alpha$ and~$c$ by $X=c$ and $Y=\alpha+c+1$.
\begin{defn}
The associated Laguerre polynomials $L_n(x;X,Y)$ are given by~\eqref{3rr} with
\begin{gather*}
b_n=2n+X+Y,
\qquad
\lambda_n=(n+X)(n+Y-1).
\end{gather*}
\end{defn}

The double sum expression for $L_n(x;X,Y)$ we shall prove is the following.
We use the notation
\begin{gather*}
(n)_k = n(n+1)(n+2)\cdots(n+k-1).
\end{gather*}

\begin{thm}
\label{doublesumthm}
For any non-negative integer~$n$,
\begin{gather*}
L_n(x;X,Y)=\sum\limits_{k=0}^n (X+k+1)_{n-k}\binom{n}{k} \, \hyper{3}{2}{k-n, k+1,Y-1}{-n,X+k+1}{1} (-1)^{n-k} x^k.
\end{gather*}
\end{thm}

We shall prove Theorem~\ref{doublesumthm} and its relative Theorem~\ref{doublesumthm2} using the even-odd polynomials
in~\cite[p.~40]{Chi}.
A~bonus of this approach is that a~closely related set of associated Laguerre polynomials, referred to as Model~II~\cite[p.~158]{Ism}
may be simultaneously studied.

\begin{defn}
The Model II associated Laguerre polynomials $L_n^{(2)}(x;X,Y)$ are given by~\eqref{3rr} with
\begin{gather*}
b_0=X,
\qquad
b_n=2n+X+Y-1,
\qquad
n>0,
\qquad
\lambda_n = (n-1+X)(n-1+Y).
\end{gather*}
\end{defn}

\begin{thm}
\label{doublesumthm2}
For any non-negative integer~$n$,
\begin{gather*}
L_n^{(2)}(x;X,Y)= \sum\limits_{k=0}^n(X+k)_{n-k}\binom{n}{k}\,  \hyper32{k-n,k,Y-1}{-n,X+k}1 (-1)^{n-k} x^k.
\end{gather*}
\end{thm}

We now review~\cite[p.~40]{Chi},~\cite[p.~311]{SimionStanton}, the basic setup for the even-odd polynomials.

Suppose that a~set of monic orthogonal polynomials $P_n(x)$ satisf\/ies~\eqref{3rr} with $b_n=0$ and $\lambda_n=\Lambda_n$
arbitrary.
Then $P_n(x)$ are even-odd polynomials, that is,
\begin{gather*}
P_{2n}(x)=e_n\big(x^2\big),
\qquad
P_{2n+1}(x)=x\thinspace o_n\big(x^2\big),
\end{gather*}
for some monic polynomials $e_n(x)$ and $o_n(x)$.

\begin{prop}
The polynomials $e_n$ and $o_n$ are both orthogonal polynomials satisfying
\begin{gather*}
e_{n+1}(x)  = (x - b_n(e)) e_n(x) - \lambda_n(e) e_{n-1}(x),
\\
o_{n+1}(x)  = (x - b_n(o)) o_n(x) - \lambda_n(o) o_{n-1}(x),
\end{gather*}
where
\begin{gather*}
b_n(e)=\Lambda_{2n+1}+\Lambda_{2n},
\qquad
\lambda_n(e)=\Lambda_{2n-1}\Lambda_{2n}
\qquad \text{and}\\
b_n(o)=\Lambda_{2n+2}+\Lambda_{2n+1},
\qquad
\lambda_n(o)=\Lambda_{2n+1}\Lambda_{2n}.
\end{gather*}
Here we define $\Lambda_0=0$, i.e.,~$b_0(e)=\Lambda_1$.
\end{prop}

We see that the choices of
%\begin{gather*}
$\Lambda_{2n}=n-1+Y$,
%\qquad
$\Lambda_{2n+1}=n+X$,
%\end{gather*}
give the associated Laguerre polynomials
\begin{gather*}
o_n(x)=L_n(x;X,Y),
\qquad
e_n(x)=L_n^{(2)}(x;X,Y).
\end{gather*}

We shall give combinatorial interpretations for $o_n(x)$ and $e_n(x)$ which will prove Theorem~\ref{doublesumthm} and
Theorem~\ref{doublesumthm2}.
Let
\begin{gather*}
e_n(x)=\sum\limits_{k=0}^n E_{n,k} x^k (-1)^{n-k},
\qquad
o_n(x)=\sum\limits_{k=0}^n O_{n,k} x^k (-1)^{n-k}.
\end{gather*}
Note that we have the following recursive relations
\begin{gather*}
%\label{maineq}
O_{n,k}=E_{n,k}+\Lambda_{2n} O_{n-1,k},
\qquad
E_{n,k}=O_{n-1,k-1}+\Lambda_{2n-1} E_{n-1,k}.
\end{gather*}

In order to give combinatorial interpretations for $O_{n,k}$ and $E_{n,k}$, we need several def\/initions.

\begin{defn}
A~\emph{$k$-marked permutation} is a~permutation~$\pi$ in which the integers in an increasing subsequence of length~$k$ are
marked.
\end{defn}

\begin{exam}
%\label{examp1}
$\pi= \mk3 421\mk6 5$ is a~2-marked permutation of length 6, where $3$ and $6$ are marked.
\end{exam}

\begin{defn}
Let $\pi=\pi_1\pi_2\cdots\pi_n$ be a~word with distinct integers.
Then $\pi_i$ is called a~\emph{right-to-left minimum} of~$\pi$ if $\pi_i<\pi_j$ for all $j=i+1,i+2,\dots,n$.
We denote by $\mathrm{RLMIN}(\pi)$ the set of right-to-left minima of~$\pi$.
Similarly, we def\/ine $\mathrm{RLMAX}(\pi)$, $\mathrm{LRMIN}(\pi)$, and $\mathrm{LRMAX}(\pi)$, which are the sets of right-to-left
maxima, left-to-right minima, and left-to-right maxima, respectively.
In other words,
\begin{gather*}
\mathrm{RLMIN}(\pi)   = \{\pi_i: \pi_i<\pi_j \;\text{for all}\; j=i+1,i+2,\dots,n\},
\\
\mathrm{RLMAX}(\pi)   = \{\pi_i: \pi_i>\pi_j \;\text{for all}\; j=i+1,i+2,\dots,n\},
\\
\mathrm{LRMIN}(\pi)   = \{\pi_i: \pi_i<\pi_j \;\text{for all}\; j=1,2,\dots,i-1\},
\\
\mathrm{LRMAX}(\pi)   = \{\pi_i: \pi_i>\pi_j \;\text{for all}\; j=1,2,\dots,i-1\}.
\end{gather*}
Let $\mathrm{RLmin}(\pi)$, $\mathrm{RLmax}(\pi)$, $\mathrm{LRmin}(\pi)$, and $\mathrm{LRmax}(\pi)$ be the cardinalities of the
sets $\mathrm{RLMIN}(\pi)$, $\mathrm{RLMAX}(\pi)$, $\mathrm{LRMIN}(\pi)$, and $\mathrm{LRMAX}(\pi)$, respectively.
\end{defn}

\begin{defn}
Let $\pi=w_0 \mk{a_1} w_1 \mk{a_2} w_2 \cdots \mk{a_k}w_k$ be a~$k$-marked permutation, where $a_1<a_2<\cdots<a_k$ are the marked
integers.
Let
\begin{gather*}
\mathrm{LRMIN}'(\pi)  = \mathrm{LRMIN}(\pi)\setminus\mathrm{LRMIN}(w_0)\setminus\{a_1,a_2,\dots,a_k\},
\\
\mathrm{RLMIN}'(\pi)  = \bigcup_{i=0}^k \mathrm{RLMIN}(w_i) \cap [a_i,a_{i+1}],
\end{gather*}
where $a_0=0$, $a_{k+1}=\infty$, and $[a_i,a_{i+1}]=\{a_i,a_i+1,\dots,a_{i+1}\}$.
We also denote $\mathrm{LRmin}'(\pi)=|\mathrm{LRMIN}'(\pi)|$ and $\mathrm{RLmin}'(\pi)=|\mathrm{RLMIN}'(\pi)|$.
\end{defn}

\begin{exam}
Let $\pi=68\mk2415\mk739$ with $k=2$, $a_1=2$, $a_2=7$, $w_0=68$, $w_1=415$, and $w_2=39$.
Then
\begin{gather*}
\mathrm{RLMIN}(\pi) =\{1,3,9\}, \qquad  \mathrm{RLMAX}(\pi) =\{9\},
\\
\mathrm{LRMIN}(\pi) =\{6,2,1\}, \qquad  \mathrm{LRMAX}(\pi) =\{6,8,9\},
\end{gather*}
and
\begin{gather*}
\mathrm{LRMIN}'(\pi) =\{6,2,1\}\setminus\{6\}\setminus\{2,7\}= \{1\},
\\
\mathrm{RLMIN}'(\pi) =\big(\{6,8\}\cap[0,2]\big)\cup \big(\{1,5\}\cap[2,7]\big) \cup\big(\{3,9\}\cap[7,\infty]\big)\\
\hphantom{\mathrm{RLMIN}'(\pi)}{} =
\varnothing\cup\{5\}\cup \{9\}=\{5,9\}.
\end{gather*}
\end{exam}

\begin{defn}
Let $\mathfrak{S}_{n,k}^E$ be the set of~$k$-marked permutations of $[n]$.
Let $\mathfrak{S}_{n,k}^O$ be the set of $(k+1)$-marked permutations of $[n+1]$ with $n+1$ marked.
\end{defn}

We now give combinatorial interpretations for the coef\/f\/icients $E_{n,k}$ and $O_{n,k}$.

\begin{prop}
\label{prop:EO}
Let $E_{n,k}$ and $O_{n,k}$ be determined by the following recursions with initial conditions $E_{n,k}=O_{n,k}=0$ if $n<k$,
$k<0$, or $n<0$; $E_{n,k}=O_{n,k}=1$ if $n=k\ge0$:
\begin{gather*}
E_{n,k}=O_{n-1,k-1}+(n-1+X)E_{n-1,k},
\qquad
O_{n,k}=E_{n,k}+(n-1+Y)O_{n-1,k}.
\end{gather*}
Then
\begin{gather*}
E_{n,k} = \sum\limits_{\pi\in\mathfrak{S}_{n,k}^E} X^{\mathrm{RLmin}'(\pi)} Y^{\mathrm{LRmin}'(\pi)},
\qquad
O_{n,k} = \sum\limits_{\pi\in\mathfrak{S}_{n,k}^O} X^{\mathrm{RLmin}'(\pi)} Y^{\mathrm{LRmin}'(\pi)}.
\end{gather*}
\end{prop}
\begin{proof}
We will show that
\begin{gather*}
A_{n,k}=\sum\limits_{\pi\in\mathfrak{S}_{n,k}^E} X^{\mathrm{RLmin}'(\pi)} Y^{\mathrm{LRmin}'(\pi)},
\qquad
B_{n,k} = \sum\limits_{\pi\in\mathfrak{S}_{n,k}^O} X^{\mathrm{RLmin}'(\pi)} Y^{\mathrm{LRmin}'(\pi)}
\end{gather*}
satisfy the recurrence relations.

Let $\pi\in\mathfrak{S}_{n,k}^E$.
If~$n$ is marked in~$\pi$, then we can consider~$\pi$ as an element of $\mathfrak{S}_{n-1,k-1}^O$.
Otherwise,~$\pi$ is obtained from an element~$\sigma$ in $\mathfrak{S}_{n-1,k}^E$ by inserting an unmarked integer~$n$ in~$n$
dif\/ferent ways.
Then we always have $\mathrm{RLmin}'(\pi)=\mathrm{RLmin}'(\sigma)$ and $\mathrm{LRmin}'(\pi)=\mathrm{LRmin}'(\sigma)$ unless we
insert~$n$ at the end in which case we have $\mathrm{RLmin}'(\pi)=\mathrm{RLmin}'(\sigma)+1$ and
$\mathrm{LRmin}'(\pi)=\mathrm{LRmin}'(\sigma)$.
This implies that $A_{n,k}=B_{n-1,k-1}+(n-1+X)A_{n-1,k}$.

Now let $\pi\in\mathfrak{S}_{n,k}^O$.
If the last element of~$\pi$ is the marked $n+1$, then by removing $n+1$ we can consider~$\pi$ as an element of
$\mathfrak{S}_{n,k}^E$ without changing $\mathrm{RLmin}'(\pi)$ and $\mathrm{LRmin}'(\pi)$.
Otherwise, let~$j$ be the last element of~$\pi$.
Since marked integers form an increasing sequence and~$n+1$ is marked, $j$~must be unmarked.
By removing~$j$ and relabeling in order-preserving way, we can consider~$\pi$ as an element~$\sigma$ in $\mathfrak{S}_{n-1,k}^O$.
Then we have $\mathrm{RLmin}'(\pi)=\mathrm{RLmin}'(\sigma)$ and $\mathrm{LRmin}'(\pi)=\mathrm{LRmin}'(\sigma)$ unless $j=1$ in
which case we have $\mathrm{RLmin}'(\pi)=\mathrm{RLmin}'(\sigma)$ and $\mathrm{LRmin}'(\pi)=\mathrm{LRmin}'(\sigma)+1$.
This implies that $B_{n,k}=A_{n,k}+(n-1+Y)B_{n-1,k}$.
\end{proof}

The next proposition constitutes proofs of Theorems~\ref{doublesumthm} and~\ref{doublesumthm2}.
\begin{prop}
%\label{mainprop}
For $n,k\ge0$, we have
\begin{gather*}
O_{n,k} = (X+k+1)_{n-k}\binom{n}{k} \, \hyper32{k-n,k+1,Y-1}{-n, X+k+1}1,
\end{gather*}
and
\begin{gather*}
E_{n,k} = (X+k)_{n-k}\binom{n}{k} \, \hyper32{k-n,k,Y-1}{-n,X+k}1.
\end{gather*}
\end{prop}
\begin{proof}
By Proposition~\ref{prop:EO} we have
\begin{gather*}
E_{n,k} = \sum\limits_{\pi\in\mathfrak{S}_{n,k}^E} X^{\mathrm{RLmin}'(\pi)} Y^{\mathrm{LRmin}'(\pi)}.
\end{gather*}

Let $\pi=w_0 \mk{a_1} w_1 \mk{a_2} w_2 \cdots \mk{a_k}w_k\in\mathfrak{S}_{n,k}^E$.
Then $\min(w_0,a_1,a_2,\dots,a_k) = m+1$ for an integer $0\le m\le n-k$.
Since $a_1<a_2<\dots<a_k$, we have two cases $m+1=a_1<\min(w_0)$ and $m+1=\min(w_0)<a_1$.

Suppose that $m+1=a_1<\min(w_0)$.
In order to construct such $\pi=w_0 \mk{a_1} w_1 \mk{a_2} w_2 \cdots \mk{a_k}w_k$ we choose $a_2,a_3,\dots,a_k$ in
$\binom{n-m-1}{k-1}$ ways, select a~permutation~$\sigma$ of $[m]$, shuf\/f\/le~$\sigma$ with $a_2,a_3,\dots,a_k$, insert $a_1=m+1$ at
the beginning, and insert the remaining $n-k-m$ integers.
Clearly, $\mathrm{LRmin}'(\pi)=\mathrm{LRmin}(\sigma)$.
In this construction, suppose that $b_1<b_2<\cdots<b_{n-k-m}$ are the remaining integers to be inserted.
We insert $b_1$ f\/irst, and then $b_2$, and so on.
Then during the~$i$th insertion process $\mathrm{RLmin}'(\pi)$ is increased by $1$ if and only if $b_i$ is inserted just before
$a_{j+1}$ for the unique integer~$j$ with $a_{j}<b_i<a_{j+1}$.
Thus inserting $b_i$ contributes a~factor of $(X+k+m+i-1)$.
Summarizing these we obtain that the sum of $X^{\mathrm{RLmin}'(\pi)} Y^{\mathrm{LRmin}'(\pi)}$ for such~$\pi$'s is equal to
\begin{gather*}
\binom{n-m-1}{k-1} \binom{m+k-1}{m} (Y)_m (X+k+m)_{n-k-m}.
\end{gather*}
Here we used the well known fact $\sum\limits_{\sigma\in\mathfrak{S}_m} Y^{\mathrm{RLmin}(\sigma)} = (Y)_m$, see~\cite[p.~125, Exercise~127]{EC1}.

Similarly, one can show that the sum of $X^{\mathrm{RLmin}'(\pi)} Y^{\mathrm{LRmin}'(\pi)}$ for $\pi=w_0 \mk{a_1}w_1\mk{a_2}w_2
\cdots \mk{a_k}w_k$ with $m+1=\min(w_0)<a_1$ is equal to
\begin{gather*}
\binom{n-m-1}{k} \binom{m+k-1}{m} (Y)_m X (X+k+m+1)_{n-k-m-1}.
\end{gather*}

Thus we have
\begin{gather*}
E_{n,k} = \sum\limits_{m=0}^{n-k} \bigg[\binom{n-m-1}{k-1} \binom{m+k-1}{m} (Y)_m (X+k+m)_{n-k-m}
\\
\phantom{E_{n,k} = }{}
+\binom{n-m-1}{k} \binom{m+k-1}{m} (Y)_m X (X+k+m+1)_{n-k-m-1} \bigg],
\end{gather*}
which is an expression for $E_{n,k}$ as sum of two terminating $\F32(1)$'s.
However a~$\F32(1)$ contiguous relation, which is equivalent to the binomial relation{\samepage
\begin{gather*}
\binom{n-m-1}{k-1}\binom{m+k-1}{m} (Y)_m + \binom{n-m-1}{k}\binom{m+k-1}{m} (Y)_m
\\
\qquad{}
-\binom{n-m}{k}\binom{m+k-2}{m-1}(k+m-1) (Y)_{m-1} =\binom{n-m}{k}\binom{m+k-1}{m} (Y-1)_{m},
\end{gather*}
gives the stated formula for $E_{n,k}$ as a~single $\F32(1)$.}

By the same argument, we obtain
\begin{gather*}
O_{n,k} = \sum\limits_{m=0}^{n-k} \bigg[\binom{n-m-1}{k-1} \binom{m+k}{m} (Y)_m (X+k+m+1)_{n-k-m}
\\
\phantom{O_{n,k} = }{}
+\binom{n-m-1}{k} \binom{m+k-1}{m} (Y)_m X (X+k+m+1)_{n-k-m-1} \bigg].
\end{gather*}

Using
\begin{gather*}
\binom{n-m-1}{k-1}\binom{m+k}{m} (Y)_m + \binom{n-m-1}{k}\binom{m+k}{m} (Y)_m
\\
\qquad
-\binom{n-m}{k}\binom{m+k-1}{m-1}(k+m) (Y)_{m-1} =\binom{n-m}{k}\binom{m+k}{m} (Y-1)_{m},
\end{gather*}
one can obtain the single sum formula for $O_{n,k}$.
\end{proof}

\begin{remark}
The $\F32$ in Theorem~\ref{doublesumthm} which is the coef\/f\/icient of $x^{k}$ is not equal to the $\F32$ given
in~\cite[equation~(5.6.18)]{Ism}.
These two $\F32$'s are related by a~$\F32(1)$ transformation.
\end{remark}

\section{Moments of associated Laguerre polynomials}%\label{sec:moments-assoc-lagu}

In this section we interpret the moments for the measures of the associated Laguerre polynomials.
This is a~routine application of Viennot's~\cite{ViennotOP} theory, details are not provided.
We need the following def\/inition of {\it{pivots}} in order to give a~combinatorial description for the moments.
These are also called {\it{strong fixed points}} in~\cite[p.~125, Exercise~128(b)]{EC1}.
\begin{defn}
Let $\pi=\pi_1\pi_2\cdots\pi_n$ be a~permutation.
If $\pi_i$ is both a~left-to-right maximum and a~right-to-left minimum, then $\pi_i$ is called a~\emph{pivot} of~$\pi$.
We denote by $\mathrm{pivot}(\pi)$ the number of pivots of~$\pi$.
Note that if $\pi_i$ is a~pivot, then $\pi_i=i$, since $\pi_i$ is greater than all of $\pi_1,\pi_2,\dots,\pi_{i-1}$ and smaller
than all of $\pi_{i+1},\pi_{i+2},\dots,\pi_n$.
\end{defn}

Suppose that the orthogonal polynomials $p_n(x)$ satisfy~\eqref{3rr}.
Then the $\nth$ moment, $\mu_n$, of a~measure for $p_n(x)$ is a~polynomial in the sequences $\{b_k\}_{k\ge0}$ and
$\{\lambda_k\}_{k\ge1}$. We denote this moment~by
\begin{gather*}
\mu_n\big(\{b_k\}_{k\ge0},\{\lambda_k\}_{k\ge1}\big)
\qquad
\text{or}
\qquad
\mu_n(b_k,\lambda_k).
\end{gather*}

The f\/irst result uses a~special choice for $b_0$.

\begin{lem}
\label{lem:XYZ}
For an integer $n\ge0$, we have
\begin{gather*}
\mu_n(b_k,(k-1+X)(k-1+Y))= \sum\limits_{\pi\in\mathfrak{S}_n} X^{\mathrm{RLmin}(\pi)} Y^{\mathrm{LRmax}(\pi)}
Z^{\mathrm{pivot}(\pi)},
\end{gather*}
where $b_k=2k-1+X+Y$ for $k>0$, and $b_0=XYZ$.
\end{lem}
\begin{proof}
This can be proved using one of the bijections between Laguerre histories and permutations, for instance, see~\cite{Si-St}.
\end{proof}

When $Z=1/Y$ in Lemma~\ref{lem:XYZ}, we obtain moments for the Model II associated Laguerre polynomials.

\begin{thm}
\label{L2m}
The moments of the Model II associated Laguerre polynomials $L_n^{(2)}(x;X,Y)$ are
\begin{gather*}
\mu_n(2k-1+X+Y^{1-\delta_{k,0}},(k-1+X)(k-1+Y))= \sum\limits_{\pi\in\mathfrak{S}_n} X^{\mathrm{RLmin}(\pi)}
Y^{\mathrm{LRmax}(\pi)-\mathrm{pivot}(\pi)}.
\end{gather*}
\end{thm}

Using the substitution $(X,Y,Z)\to(X+1,Y,\frac{X+Y}{(X+1)Y})$ in Lemma~\ref{lem:XYZ}, we obtain
\begin{thm}
\label{L1m}
The moments of the Model I associated Laguerre polynomials $L_n(x;X,Y)$ are
\begin{gather*}
\mu_n(2k+X+Y,(k+X)(k-1+Y))
\\
\qquad
= \sum\limits_{\pi\in\mathfrak{S}_n}
(X+1)^{\mathrm{RLmin}(\pi)-\mathrm{pivot}(\pi)}(X+Y)^{\mathrm{pivot}(\pi)} Y^{\mathrm{LRmax}(\pi)-\mathrm{pivot}(\pi)}.
\end{gather*}
\end{thm}

If $X=Y=1$ in Theorem~\ref{L1m}, we have
\begin{gather*}
\mu_n(2k+2,k(k+1))=(n+1)!.
\end{gather*}
It is therefore more natural to write the $\nth$ moment of the Model~I associated Laguerre polynomial as a~sum over
$\mathfrak{S}_{n+1}$ instead of $\mathfrak{S}_{n}$. The next lemma, which is a~restatement of the even-odd polynomial
construction in Section~\ref{sec:assoc-lagu-polyn}, accomplishes this.

\begin{lem}
\label{lem:evenodd}
Given a~sequence $\{a_k\}_{k\ge0}$ with $a_0=0$, we have
\begin{gather*}
\mu_{2n}(0, a_k) = \mu_n(a_{2k}+a_{2k+1}, a_{2k-1}a_{2k}),
\\
\mu_{2n+2}(0, a_k) = a_1 \mu_n(a_{2k+1}+a_{2k+2}, a_{2k}a_{2k+1}).
\end{gather*}
\end{lem}

\begin{thm}
\label{thm:Lmu}
The moments of the Model I associated Laguerre polynomials $L_n(x;X,Y)$ are
\begin{gather*}
\mu_n(2k+X+Y,(k+X)(k-1+Y))= \sum\limits_{\pi\in\mathfrak{S}_{n+1}} X^{\mathrm{RLmin}(\pi)-1}
Y^{\mathrm{LRmax}(\pi)-\mathrm{pivot}(\pi)}.
\end{gather*}
\end{thm}
\begin{proof}
Let $a_{2k-1}=k-1+X$ and $a_{2k}=k-1+Y$ for $k\ge1$ and $a_0=0$.
Then by Lemma~\ref{lem:evenodd}, we have
\begin{gather*}
\mu_{2n}(0,a_k) = \mu_n\big(2k-1+X+Y^{1-\delta_{k,0}},(k-1+X)(k-1+Y)\big),
\\
\mu_{2n+2}(0,a_k)=X\mu_n(2k+X+Y,(k+X)(k-1+Y)).
\end{gather*}
Thus, we have
\begin{gather*}
X\mu_n(2k+X+Y,(k+X)(k-1+Y))
\\
\qquad
= \mu_{n+1}\big(2k-1+X+Y^{1-\delta_{k,0}},(k-1+X)(k-1+Y)\big).
\end{gather*}
By applying Theorem~\ref{L2m} to the above equation, we obtain the desired identity.
\end{proof}

As a~corollary of the proof above, we obtain the following generating function identity.

\begin{cor}
For any non-negative integer~$n$, we have
\begin{gather*}
\sum\limits_{\pi\in\mathfrak{S}_n} (X+1)^{\mathrm{RLmin}(\pi)-\mathrm{pivot}(\pi)}(X+Y)^{\mathrm{pivot}(\pi)}
Y^{\mathrm{LRmax}(\pi)-\mathrm{pivot}(\pi)}
\\
\qquad
=\sum\limits_{\pi\in\mathfrak{S}_{n+1}} X^{\mathrm{RLmin}(\pi)-1} Y^{\mathrm{LRmax}(\pi)-\mathrm{pivot}(\pi)}.
\end{gather*}
\end{cor}

We may analogously push the $n!$ terms for the Model II moments to $\mathfrak{S}_{n-1}$.
\begin{cor}
The moments of the Model II associated Laguerre polynomials $L_n^{(2)}(x;X,Y)$ are
\begin{gather*}
\mu_n(b_k,(k-1+X)(k-1+Y))
\\
\qquad
=(X+1)\sum\limits_{\pi\in\mathfrak{S}_{n-1}} (X+2)^{\mathrm{RLmin}(\pi)-\mathrm{pivot}(\pi)}(X+Y+1)^{\mathrm{pivot}(\pi)}
Y^{\mathrm{LRmax}(\pi)-\mathrm{pivot}(\pi)},
\end{gather*}
where $b_k=2k-1+X+Y$ for $k>0$, and $b_0=X$.
\end{cor}

\section{Permutation tableaux}%\label{sec:permutation-tableaux}

In this section we reinterpret the weights on permutations of the previous section to be weights on another set of objects also
counted by $n!$, called permutation tableaux.
Permutation tableaux were introduced by Postnikov~\cite{Postnikov} in the study of totally nonnegative Grassmanian and studied
extensively.

A \emph{permutation tableau} is a~f\/illing of a~Ferrers diagram using 0's and 1's such that each column contains at least one 1
and there is no 0 which has a~1 above in the same column and a~1 to the left of it in the same row.
We allow a~Ferrers diagram to have empty rows.

\begin{figure}
\centering \includegraphics{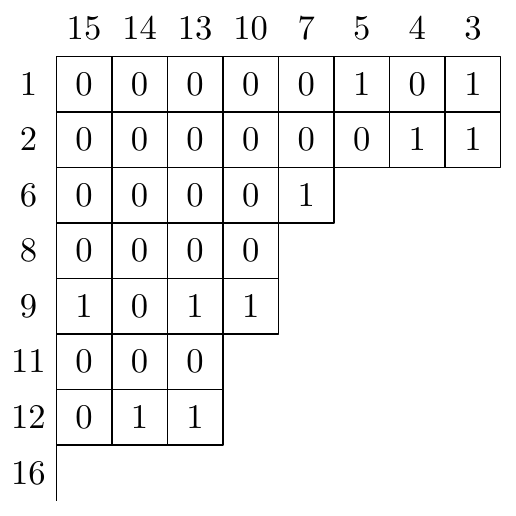}
\qquad
\includegraphics{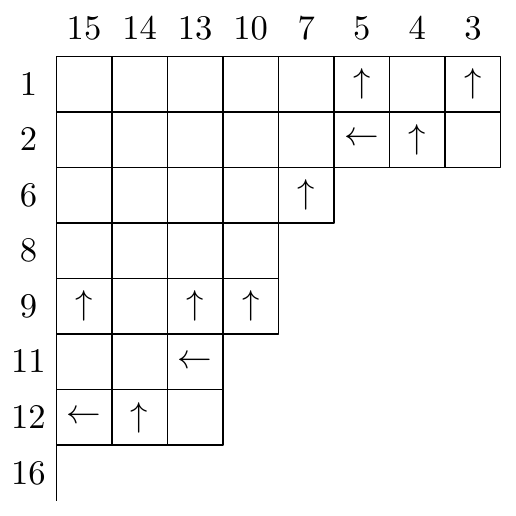} \caption{A permutation tableau on the left and its alternative representation on the right.}
\label{fig:pt}
\end{figure}

\looseness=-1
The \emph{length} of a~permutation tableau is the number of rows plus the number of columns.
Let $\mathrm{PT}_n$ denote the set of permutation tableau of length~$n$.
If a~permutation tableau has length~$n$, we label the southeast border of a~Ferrers diagram with $1,2,\dots,n$ starting from
northeast to southwest.
The rows and columns are also labeled by the same labels as in the border steps contained in the rows and columns.
An example of a~permutation tableau is shown on the left of Fig.~\ref{fig:pt}.

In a~permutation tableau, a~\emph{topmost $1$} is a~1 which is the topmost 1 in the column containing it.
A~\emph{restricted~$0$} is a~0 which has a~1 above it in the same column.
A~\emph{rightmost restricted~$0$} is the rightmost restricted~0 in the row containing it.
A~row is called \emph{restricted} if it has a~restricted~0, and \emph{unrestricted} otherwise.
We denote by $\mathrm{urr}(T)$ the number of unrestricted rows in~$T$.
Similarly, a~0 in $T\in\mathrm{PT}_n$ is called \emph{c-restricted} if it has a~1 to the left of it in the same row.
A~column is called \emph{restricted} if it has a~c-restricted~0, and \emph{unrestricted} otherwise.
We denote by~$\mathrm{urc}(T)$ the number of unrestricted columns in~$T$.
For example, if~$T$ is the permutation tableau in Fig.~\ref{fig:pt}, then the unrestricted rows are rows~1,~6,~8,~9,~16, and the
unrestricted columns are columns~3,~5,~7,~10,~13,~15.

Corteel and Nadeau~\cite{Corteel2009} found a~bijection $\Phi\colon \mathrm{PT}_n\to\mathfrak{S}_n$ such that if $\Phi(T)=\pi$, then
$\mathrm{urr}(T)=\mathrm{RLmin}(\pi)$.
Corteel and Kim~\cite{CorteelKim} found a~generating function for permutation tableaux according to the number of unrestricted
columns.

The \emph{alternative representation} of a~permutation tableau is the f\/illing of the same Ferrers diagram obtained by replacing
each topmost~1 by an up arrow, each rightmost restricted 0 by a~left arrow, and leaving the other entries empty.
See Fig.~\ref{fig:pt}.
The diagram obtained in this way is called an \emph{alternative tableau}.
Equivalently, an alternative tableau can be def\/ined as a~f\/illing of a~Ferrers diagram with uparrows and leftarrows such that each
column contains exactly one uparrow and there is no arrow pointing to another arrow.
Alternative tableaux were f\/irst considered by Viennot~\cite{Viennota} and studied further by Nadeau~\cite{Nadeau}.

In this section we prove the following theorem.

\begin{thm}
\label{PT}
The moments for the associated Laguerre polynomials are
\begin{gather*}
\mu_n(2k+X+Y,(k+X)(k-1+Y))=\sum\limits_{T\in\mathrm{PT}_{n+1}} X^{\mathrm{urr}(T)-1} Y^{\mathrm{urc}(T)}.
\end{gather*}
\end{thm}

For the remainder of this section we shall prove the above theorem.

It is easy to see that in the alternative representation of $T\in\mathrm{PT}_n$, a~row is unrestricted if and only if it has no
left arrow.
There is also a~simple characterization for unrestricted columns in the alternative representation.

\begin{lem}
\label{lem:urc}
In the alternative representation of $T\in\mathrm{PT}_n$, a~column is unrestricted if and only if there is no arrow that is
strictly to the northwest of the uparrow in the column.
\end{lem}
\begin{proof}
Suppose that a~column, say column~$c$, is restricted.
Then column~$c$ has a~c-restric\-ted~0.
Suppose that this c-restricted 0 has a~1 in column $c'$ to the left of it in the same row.
Then by the condition for permutation tableaux, column~$c$ cannot have a~1 above the c-restricted~0.
Thus the topmost 1 in column $c'$ is strictly to the northwest of the topmost~1 in column~$c$.
Since topmost~1s correspond to uparrows in the alternative representation, the uparrow in column~$c$ has an arrow strictly to the
northwest of it.
The converse can be proved similarly.
\end{proof}

For a~permutation $\pi=\pi_1\cdots\pi_n$, we def\/ine $\pi^{\mathrm{rev}}=\pi_n\pi_{n-1}\cdots\pi_1$.
We also def\/ine $\pi^*$ as follows.
If $\pi=w_1 r_1 w_2 r_2 \cdots w_kr_k$, where $r_1,r_2,\dots,r_k$ are the right-to-left minima of~$\pi$, then $\pi^* =
w_1^{\mathrm{rev}} r_1 w_2^{\mathrm{rev}} r_2 \cdots w_k^{\mathrm{rev}} r_k$.
Note that $\pi^*$ and~$\pi$ have the same right-to-left minima and \mbox{$(\pi^*)^*=\pi$}.
For example, if
\begin{gather*}
%\label{eq:pi}
\pi = 4,2,5,3,\overline{1},7,\overline{6},\overline{8},14,12,15,11,13,10,\overline{9},\overline{16},
\end{gather*}
where the right-to-left minima are overlined, then
\begin{gather*}
%\label{eq:pi*}
\pi^*=\underline{3},\underline{5},2,4,\overline{1},\underline{7},
\overline{6},\underline{\overline{8}},\underline{10},\underline{13},11,\underline{15},
12,14,\overline{9},\underline{\overline{16}},
\end{gather*}
where the right-to-left minima are overlined and the left-to-right maxima are underlined.

We now describe the bijection $\Phi\colon \mathrm{PT}_n\to\mathfrak{S}_n$ of Corteel and Nadeau~\cite{Corteel2009} using the
alternative representation.
Let $T\in\mathrm{PT}_n$.
Suppose that $r_1<r_2<\dots<r_k$ are the labels of the unrestricted rows of~$T$.
We f\/irst set~$\pi$ to be $r_1r_2\cdots r_k$.
For each column of~$T$ starting from left to right we will insert integers to~$\pi$ as follows.
If column~$i$ contains the uparrow in row~$j$ and the leftarrows in rows $i_1<\dots<i_r$, then insert $i_1,\dots,i_r,i$ in this
order before~$j$ in~$\pi$.
If column~$i$ does not have any leftarrows, we insert only~$i$ before~$j$.
Observe that we have $j<i_1<i_2<\dots<i_r<i$.
Thus inserting integers in this way does not change the right-to-left minima.
We def\/ine $\Phi(T)$ to be the resulting permutation~$\pi$.

We will use the following property of the map~$\Phi$.
If $\Phi(T)=\pi$, then~$i$ is the label of a~column if and only~$i$ is followed by a~smaller integer in~$\pi$.

\begin{lem}
Let $\Phi(T)=\pi$.
Then we have the following.
\begin{itemize}\itemsep=0pt
\item An integer~$r$ is the label of an unrestricted row of~$T$ if and only if it is a~right-to-left minimum of $\pi^*$.
\item An integer~$c$ is the label of an unrestricted column of~$T$ if and only if it is a~left-to-right maximum but not a~pivot
of $\pi^*$.
\end{itemize}
\end{lem}
\begin{proof}
From the construction of the bijection~$\Phi$, the labels of unrestricted rows of~$T$ are the right-to-left minima of~$\pi$.
Since~$\pi$ and $\pi^*$ have the same right-to-left minima, we obtain the f\/irst item.

We now prove the second item.
Suppose that column~$c$ is an unrestricted column of~$T$.
Let row~$r$ be the row containing the uparrow of column~$c$.
By Lemma~\ref{lem:urc}, row~$r$ must be unrestricted.
Let $r_1<r_2<\cdots<r_k$ be the labels of unrestricted rows of~$T$, and let $r=r_i$.
By the construction of the map~$\Phi$, before inserting integers from column~$c$, the integers that have been inserted are
located to the right of $r_{i-1}$.
Then after inserting integers from column~$c$, we have~$c$ to the left of $r_i$ with no integers between them.
Since the remaining integers to be inserted are less than~$c$, every integer greater than~$c$ is located either to the right of~$r_i$ or between~$r_{i-1}$ and~$c$ in $\pi=\Phi(T)$.
Thus~$c$ is a~left-to-right maximum of~$\pi^*$.
Since~$c$ is to the left of~$r_i$ and $c<r_i$,~$c$ is not a~right-to-left minimum of~$\pi^*$.

Now let~$c$ be a~left-to-right maximum of $\pi^*$ that is not a~right-to-left minimum.
Let $r_1<r_2<\dots<r_k$ be the right-to-left minima of~$\pi$ (equivalently~$\pi^*$).
Suppose that~$c$ is located between~$r_{i-1}$ and~$r_i$.
It is easy to see that~$c$ is followed by a~smaller number in~$\pi$.
Thus~$c$ is the label of a~column in~$T$.
If~$c$ is a~restricted column, then by Lemma~\ref{lem:urc}, there is an arrow strictly to the northwest of the uparrow in
column~$c$.
Then the column containing such an arrow has label $c'>c$.
When we insert integers from column $c'$ in the construction of $\Phi(T)$, $c'$ is located to the left of $r_{i-1}$.
Then~$c$ cannot be a~left-to-right maximum of $\pi^*$, which is a~contradiction.
Thus column~$c$ is unrestricted.
\end{proof}

The above lemma implies the following.

\begin{prop}
\label{PTprop}
We have
\begin{gather*}
\sum\limits_{T\in\mathrm{PT}_n} X^{\mathrm{urr}(T)} Y^{\mathrm{urc}(T)} =\sum\limits_{\pi\in\mathfrak{S}_n}
X^{\mathrm{RLmin}(\pi)} Y^{\mathrm{LRmax}(\pi)-\mathrm{pivot}(\pi)}.
\end{gather*}
\end{prop}

Theorem~\ref{PT} follows from Theorem~\ref{thm:Lmu} and Proposition~\ref{PTprop}.

\section{Proof via the moment generating function}\label{sec:proof-via-moment}

In this section we prove the following equivalent version of Theorem~\ref{doublesumthm}
\begin{gather}
\label{eq:3F2}
L_n(x;X,Y)=\frac{(Y)_n(X+1)_n(-1)^n}{n!} \sum\limits_{k=0}^n \frac{(-n)_k}{(Y)_k(X+1)_k} \, \hyper32{k-n,X,Y-1}{Y+k,X+k+1}1 x^k.
\end{gather}

We use only specialized knowledge of the moments.
What we need is a~recurrence relation for the moments.
For ease of notation let's put
\begin{gather*}
\theta_n(X,Y)= \mu_n(2k+X+Y,(k+X)(k-1+Y)).
\end{gather*}

\begin{prop}
The moments of the associated Laguerre polynomials satisfy the recurrence relation
\begin{gather}
\label{murelation}
\sum\limits_{k=0}^J \theta_k(X,Y) \frac{(Y-1)_{J-k}(X)_{J-k}}{(J-k)!}= \frac{(Y)_J (X+1)_J}{J!},
\qquad
J\ge 0.
\end{gather}
\end{prop}

\begin{proof}
The recurrence~\eqref{murelation} is equivalent to the explicit form of the moment generating function, which is
Proposition~\ref{contfrac}, a~quotient of hypergeometric series.
This appears in~\cite[Theo\-rem~6.5, p.~212]{JT},
or may easily be verif\/ied using contiguous relations for $\F20(x)$ to f\/ind the relevant Jacobi continued fraction.
\end{proof}

\begin{prop}
\label{contfrac}
The generating function, as a~formal power series in~$x$, for the associated Laguerre moments is
\begin{gather*}
\sum\limits_{n=0}^\infty \mu_n(2k+X+Y,(k+X)(k-1+Y)) x^n= \frac{\F20(Y,X+1;x)}{\F20(Y-1,X;x)}.
\end{gather*}
\end{prop}

Let $\mathcal{L}$ be the linear functional which uses the moments $\theta_n(X,Y)$, i.e., $\mathcal{L}(x^n)=\theta_n(X,Y)$.
We prove~\eqref{eq:3F2} by showing that the linear functional $\mathcal{L}$ satisf\/ies
\begin{gather}
\label{orth}
\mathcal{L}(x^s L_n(x;X,Y))=0,
\qquad
\text{for}
\qquad
0\le s\le n-1.
\end{gather}
We can rewrite $L_n(x;X,Y)$ as
\begin{gather*}
L_n(x;X,Y)= \frac{(Y)_n (X+1)_n (-1)^n}{n!}\sum\limits_{J=0}^n \frac{(-n)_J}{(Y)_J (X+1)_J} \sum\limits_{k=0}^J \frac{(Y-1)_{J-k}
(X)_{J-k}}{(J-k)!} x^k.
\end{gather*}

First we establish the $s=0$ case of~\eqref{orth}.
Applying the linear functional $\mathcal{L}$ and using~\eqref{murelation}, we obtain, for $n\ge1$,
\begin{gather*}
\mathcal{L}(L_n(x;X,Y)) = \frac{(Y)_n (X+1)_n (-1)^n}{n!}\sum\limits_{J=0}^n \frac{(-n)_J}{(Y)_J (X+1)_J} \sum\limits_{k=0}^J
\frac{(Y-1)_{J-k} (X)_{J-k}}{(J-k)!} \theta_k(X,Y)
\\
\phantom{\mathcal{L}(L_n(x;X,Y))}
 =\frac{(Y)_n (X+1)_n (-1)^n}{n!}\sum\limits_{J=0}^n \frac{(-n)_J}{(Y)_J (X+1)_J} \frac{(Y)_J (X+1)_J}{J!}
\\
\phantom{\mathcal{L}(L_n(x;X,Y))}
 = \frac{(Y)_n (X+1)_n (-1)^n}{n!}\sum\limits_{J=0}^n \binom{n}{J} (-1)^J=0.
\end{gather*}

We'll f\/inish proving~\eqref{orth} for $0< s<n$ by obtaining, at the last step, an $\nth$ dif\/ference of a~polynomial of degree
$s$, which is zero.

Replacing $\theta_k(X,Y)$ in the above computation with $\theta_{k+s}(X,Y)$ and reversing the interior sum, shows that we need
the following lemma.

\begin{lem}
For $s\ge 1$ an integer,
\begin{gather*}
\sum\limits_{k=0}^J \theta_{J+s-k}(X,Y) \frac{(Y-1)_k (X)_k}{k!}= \left({\text{a polynomial in~$J$ of degree~$s$}}\right)
\frac{(Y)_J (X+1)_J}{J!}.
\end{gather*}
\end{lem}

\begin{proof}
Note that~\eqref{murelation} implies that
\begin{gather*}
\sum\limits_{k=0}^J \theta_{J+s-k}(X,Y) \frac{(Y-1)_k (X)_k}{k!}= \frac{(Y)_J(X+1)_J}{J!(J+1)\cdots (J+s)} \biggl((Y+J)_s(X+J+1)_s
\\
\qquad
-\sum\limits_{p=0}^{s-1} \theta_p(X,Y) (J+s)\cdots (J+s-p+1) (Y-1)(Y+J)_{s-p-1}X(X+J+1)_{s-p-1}\biggr).
\end{gather*}

We must show that the polynomial $p(J)$ in~$J$ of degree $2s$
\begin{gather*}
p(J)= (Y+J)_s(X+J+1)_s
\\
\phantom{p(J)=}{}
-\sum\limits_{p=0}^{s-1} \theta_p(X,Y) (J+s)\cdots (J+s-p+1) (Y-1)(Y+J)_{s-p-1}X(X+J+1)_{s-p-1}
\end{gather*}
is divisible by the~$s$ factors $J+1,\dots,J+s$. To this end, put $J=-k$, $1\le k\le s$, and evaluate
\begin{gather*}
p(-k)=(Y-k)_k (X+1-k)_k\biggl((Y)_{s-k}(X+1)_{s-k}
\\
\phantom{p(-k)=}{}
-\sum\limits_{p=0}^{s-k} \theta_p(X,Y)\frac{(s-k)!}{(s-k-p)!}(Y-1)_{s-p-k}(X)_{s-p-k}\biggr) = 0
\end{gather*}
by~\eqref{murelation}.
\end{proof}

\section{Remarks}

\vspace{-1.5mm}

The combinatorics of the associated classical polynomials was f\/irst studied
by Drake who considers the associated Hermite
polynomials in~\cite{Drake2009}.
In this paper we study the combinatorics of the associated Laguerre polynomials.
Wimp~\cite[Theorem 1]{Wimp} gives an analogous double sum formula for the associated Jacobi polynomials using a~${}_4F_3$. Since
the normalized three term recurrence relation contains rational functions instead of polynomials, the combinatorics of the
moments will be more dif\/f\/icult in the Jacobi case.

In the~$q$-world, one would like to give an alternative proof of the double sum formula for the associated Askey--Wilson
polynomials involving a~${}_{10}\phi_9$, see Ismail and Rahman~\cite[equation~(4.15)]{IR}.
One could hope for a~proof along the lines of Section~\ref{sec:proof-via-moment}.
Perhaps a~recurrence for the generalized moments of the Askey--Wilson basis of type~\eqref{murelation} exists.

A~$q$-analogue of Theorem~\ref{PT} has been given by Corteel and Josuat-Verg\`es~\cite{CJV}.
One may hope that it would give intuition for the appropriate parametrization for the associated Askey--Wilson polynomials and
ASEP, see~\cite{CW}.

\vspace{-2.5mm}

\subsection*{Acknowledgements}

\vspace{-1.5mm}

The f\/irst author was partially supported by Basic Science Research Program through the National Research Foundation of Korea
(NRF) funded by the Ministry of Education (NRF-2013R1A1A2061006).
The second author was supported by NSF grant DMS-1148634.

\vspace{-2.5mm}

\pdfbookmark[1]{References}{ref}
\LastPageEnding

\end{document}